\documentclass[12pt,vatola]{article}

\usepackage{amsfonts}

\usepackage{mathrsfs}

\usepackage{graphicx}

\def\c{\centerline}

\def\re#1{\par\hangindent\parindent\indent\llap{#1\enspace}\ignorespaces}

\def\no{\noindent}

\begin{document}

\c{\bf\large Combinatorial Structure of Manifolds with}\vskip 3mm

\c{\bf\large Poincar\'{e} Conjecture} \vskip 8mm

\c{Linfan MAO}\vskip 5mm

\c{\scriptsize (Chinese Academy of Mathematics and System Science,
Beijing 100080, P.R.China)}

\c{\scriptsize E-mail: maolinfan@163.com}

\vskip 12mm

\begin{minipage}{130mm}

\no{\bf Abstract}: {\small A manifold $M^n$ inherits a labeled
$n$-dimensional graph $\widetilde{M}[G^L]$ structure consisting of
its charts. This structure enables one to characterize fundamental
groups of manifolds, classify those of locally compact manifolds
with finite non-homotopic loops by that of labeled graphs $G^L$. As
a by-product, this approach also concludes that {\it every homotopy
$n$-sphere is homeomorphic to the sphere $S^n$ for an integer $n\geq
1$}, particularly, the Perelman's result for $n=3$.}\vskip 2mm

\no{\bf Key Words}: {\small Combinatorial Euclidean space, manifold,
fundamental group, labeled graph, homotopic loop.}\vskip 2mm

\no{\bf AMS(2010)}: {\small 51H20.}

\end{minipage}

\vskip 12mm

\no{\bf \S $1.$ \ Introduction}

\vskip 4mm

\no An $n$-manifold is a second countable Hausdorff space of locally
Euclidean $n$-space without boundary. Terminologies and notions used
in this paper are standard. We follow references [1], [13] for
topology, [5], [6] for graphs or topological graphs and [11] for
combinatorial manifolds. Here, we only mention conceptions appeared
in references [8], [10]-[11] for combinatorial manifolds.

\vskip 4mm

\no{\bf Definition $1.1$} \ {\it A combinatorial Euclidean space
$\mathscr{E}_G(n_{\nu};\nu\in\Lambda)$ underlying a connected graph
$G$ is a topological spaces consisting of ${\bf R}^{n_{\nu}}$,
$\nu\in\Lambda$ for an index set $\Lambda$ such that\vskip 3mm

$V(G)=\{{\bf R}^{n_{\nu}}|\nu\in\Lambda\}$;\vskip 2mm

$E(G)=\{\ ({\bf R}^{n_{\mu}},{\bf R}^{n_{\nu}}) | \ {\bf
R}^{n_{\mu}}\cap{\bf R}^{n_{\nu}}\not=\emptyset,
\mu,\nu\in\Lambda\}$.\vskip 2mm

If $|\Lambda|=1$, i.e., the dimension of Euclidean spaces in
$\mathscr{E}_G(n_{\nu};\nu\in\Lambda)$ are all in the same $n$, the
notation $\mathscr{E}_G(n_{\nu};\nu\in\Lambda)$ is abbreviated to
$\mathscr{E}_G(n,\ \nu\in\Lambda)$ for simplicity. }\vskip 3mm

Notice that a Euclidean space ${\bf R}^n$ is an $n$-dimensional
vector space with a normal basis
$\overline{\epsilon}_1=(1,0,\cdots,0)$,
$\overline{\epsilon}_2=(0,1,0\cdots,0)$, $\cdots$,
$\overline{\epsilon}_n=(0,\cdots,0,1)$, namely, it has $n$
orthogonal orientations. In Definition $1.1$, we do not assume ${\bf
R}^{n_{\mu}}\cap {\bf R}^{n_{\nu}}={\bf
R}^{\min\{n_{\mu},n_{\nu}\}}$. In fact, let $\mathcal{X}_{v_{\mu}}$
be the set of orthogonal orientations in ${\bf R}^{n_{v_{\mu}}}$,
$\mu\in\Lambda$, respectively ([12]). Then

$${\bf R}^{n_{\mu}}\bigcap {\bf R}^{n_{\nu}}=\mathcal{X}_{v_{\mu}}\bigcap\mathcal{X}_{v_{\nu}}.$$

 A {\it combinatorial fan-space}
$\widetilde{\bf R}(n_{\nu};\nu\in\Lambda)$ is a combinatorial
Euclidean space $\mathscr{E}_{K_{|\Lambda|}}(n_{\nu};\nu\in\Lambda)$
of ${\bf R}^{n_{\nu}},\ \nu\in\Lambda$ such that for any integers
$\mu,\nu\in\Lambda,\ \mu\not= \nu$,

$${\bf R}^{n_{\mu}}\bigcap{\bf R}^{n_{\nu}}= \bigcap\limits_{\lambda\in\Lambda}{\bf
R}^{n_{\lambda}}.$$

If $|\Lambda|=m <+\infty$, for $\forall p\in\widetilde{{\bf
R}}(n_{\nu}; \nu\in\Lambda)$ we can present it by an $m\times n_{m}$
coordinate matrix $[\overline{x}]$ following with
$x_{il}=\frac{x_{l}}{m}$ for $1\leq i\leq m, 1\leq
l\leq\widehat{m}$, where $\widehat{m}$ is the number of orthogonal
orientations in
$\mathcal{X}_{v_{\mu}}\bigcap\mathcal{X}_{v_{\nu}}$,\vskip 2mm

\[
[\overline{x}]=\left[
\begin{array}{cccccccc}
x_{11} & \cdots & x_{1\widehat{m}}
& x_{1(\widehat{m})+1)} & \cdots & x_{1n_1} & \cdots & 0 \\
x_{21} & \cdots & x_{2\widehat{m}}
& x_{2(\widehat{m}+1)} & \cdots & x_{2n_2} & \cdots & 0  \\
\cdots & \cdots & \cdots  & \cdots & \cdots & \cdots  \\
x_{m1} & \cdots & x_{m\widehat{m}} & x_{m(\widehat{m}+1)} & \cdots &
\cdots & x_{mn_{m}-1} & x_{mn_{m}}
\end{array}
\right],
\]\vskip 2mm

\no which enables us to generalize the conception of manifold to
combinatorial manifold, a locally combinatorial Euclidean space.
\vskip 4mm

\no{\bf Definition $1.2$} \ {\it For a given integer sequence $0<
n_1<n_2<\cdots< n_m$, $m\geq 1$, a topological combinatorial
manifold $\widetilde{M}$ is a {\it Hausdorff space} such that for
any point $p\in \widetilde{M}$, there is a local chart
$(U_p,\varphi_p)$ of $p$, i.e., an open neighborhood $U_p$ of $p$ in
$\widetilde{M}$ and a homoeomorphism $\varphi_p:
U_p\rightarrow\widetilde{\bf R}(n_1(p),
n_2(p),\cdots,n_{s(p)}(p))=\bigcup\limits_{i=1}^{s(p)}{\bf
R}^{n_i(p)}$ with\vskip 3mm

\hskip 5mm$\{n_1(p),
n_2(p),\cdots,n_{s(p)}(p)\}\subseteq\{n_1,n_2,\cdots,n_m\}$
and\vskip 3mm

\hskip 5mm$\bigcup\limits_{p\in\widetilde{M}}\{n_1(p),
n_2(p),\cdots,n_{s(p)}(p)\}=\{n_1,n_2,\cdots,n_m\}$,\vskip 2mm

\no denoted by $\widetilde{M}(n_1,n_2,\cdots,n_m)$ or
$\widetilde{M}$ on the context and\vskip 3mm

\c{$\widetilde{{\mathcal A}}=\{(U_p,\varphi_p)|
p\in\widetilde{M}(n_1,n_2,\cdots,n_m))\}$}\vskip 2mm

\no an atlas on $\widetilde{M}(n_1,n_2,\cdots,n_m)$.

A topological combinatorial manifold
$\widetilde{M}(n_1,n_2,\cdots,n_m)$ is finite if it is just combined
by finite manifolds without one manifold contained in the union of
others.}\vskip 3mm

If $m=1$, then $\widetilde{M}(n_1,n_2,\cdots,n_m)$ is exactly the
manifold $M^{n_1}$ by definition. Furthermore, if these manifolds
$M_i, \ 1\leq i\leq m$ in $\widetilde{M}(n_1,n_2,\cdots,n_m)$ are
Euclidean spaces ${\bf R}^{n_i}, \ 1\leq i\leq m$, then
$\widetilde{M}(n_1,n_2,\cdots,n_m)$ is nothing but the combinatorial
Euclidean space $\mathscr{E}_{G}(n_{\nu};\nu\in\Lambda)$ with
$\Lambda=\{1,2,\cdots,m\}$.

For a finitely combinatorial manifold
$\widetilde{M}(n_1,n_2,\cdots,n_m)$ consisting of manifolds $M_i,\
1\leq i\leq m$, we can construct a vertex-edge labeled graph
$G^L[\widetilde{M}]$ defined by

\vskip 3mm

\hskip 5mm $V(G^L[\widetilde{M}])=\{M_1,M_2,\cdots,M_m\},$\vskip 2mm

\hskip 5mm $E(G^L[\widetilde{M})=\{\ (M_i,M_j)\ | \ M_i\bigcap
M_j\not=\emptyset, 1\leq i,j\leq n\}$\vskip 2mm

\no with a labeling mapping $\Theta: V(G^L[\widetilde{M}])\cup
E(G^L[\widetilde{M}])\rightarrow{\bf Z}^+$ determined by

$$\Theta(M_i)={\rm dim}M_i \ \ {\rm and} \ \ \Theta(M_i,M_j)={\rm
dim}M_i\bigcap M_j$$

\no for integers $1\leq i,j\leq m$, which is inherent structure of
combinatorial manifolds. Particularly, for a combinatorial Euclidean
space $\mathscr{E}_G(n_1,\cdots,n_m)$, this vertex-edge labeled
graph is defined on $G$ by labeling $\Theta: V(G)\cup
E(G)\rightarrow{\bf Z}^+$ determined by

$$\Theta({\bf R}^{n_i})=n_i \ \ {\rm and} \ \ \Theta({\bf R}^{n_i},{\bf R}^{n_j})={\rm
dim}{\bf R}^{n_i}\bigcap {\bf R}^{n_j}$$

\no for integers $1\leq i,j\leq m$.

The objective of this paper is to characterize the inherent
combinatorial structure of $n$-manifolds with finite non-homotopic
loops. For such manifolds, there is a well-known Poincar\'{e}
conjecture first proved by Perelman ([3]) following.

\vskip 4mm

\no{\bf Theorem $1.3$}(Perelman,[15]-[17]) \ {\it Any closed simply
connected $3$-manifold is homeomorphic to $S^3$.}\vskip 3mm

Notice that a {\it homotopy sphere} is an $n$-manifold homotopy
equivalent to the $n$-sphere. A generalized Poincar\'{e} conjecture
says that {\it any homotopy $n$-sphere is homeomorphic to $S^n$}.
Combining works of Smale in 1961 for $n\geq 5$ ([18]), Freedman's in
1982 for $n=4$ ([2]), Theorem $1.3$ and classical result for
$n=1,2$, we conclude that

\vskip 4mm

\no{\bf Theorem $1.4$} \ {\it Any homotopy $n$-sphere is
homeomorphic to $S^n$ for $n\geq 1$.}\vskip 3mm

Then {\it can we find a unified proof on Theorem $1.4$}? By applying
a combinatorial notion [9], we find an inherent labeled graph
$G^L[M^n]$ structure consisting of charts homeomorphic to a
Euclidean space ${\bf R}^n$ for a topological manifold $M^n$ in this
paper, which enables us to classify locally compact $n$-manifolds
with finite non-homotopic loops by labeled graphs for integers
$n\geq 1$, also simplifies the calculation of fundamental groups of
$n$-manifolds and get a combinatorial proof for Theorem $1.4$.

\vskip 12mm

\no{\bf\S$2.$ \ Dimensional Graphs}\vskip 5mm

\no We discuss combinatorial Euclidean spaces $\mathscr{E}_G(n)$ in
this section whose importance is shown in the next result.

\vskip 4mm

\no{\bf Theorem $2.1$} \ {\it A locally compact $n$-manifold $M^n$
is a combinatorial manifold $\widetilde{M}_{G}(n)$ homeomorphic to a
Euclidean space $\mathscr{E}_{G'}(n)$ with countable graphs $G\cong
G'$ inherent in $M^n$, denoted by $G[M^n]$.}

\vskip 3mm

{\it Proof} \ Let $M^n$ be a  locally compact $n$-manifold with an
atlas

$$\mathscr{A}[M^n]=\{\ (U_{\lambda};\varphi_{\lambda})\ | \
\lambda\in\Lambda\},$$

\no where $\Lambda$ is a countable set. Then each $U_{\lambda},\
\lambda\in\Lambda$ is itself an $n$-manifold by definition. Define
an underlying combinatorial structure $G$ by\vskip 3mm

\hskip 5mm $V(G)= \{U_{\lambda}| \lambda\in\Lambda\},$\vskip 2mm

\hskip 5mm $E(G)=\{\ (U_{\lambda},U_{\iota})_i,1\leq
i\leq\kappa_{\lambda\iota}+1 | \ U_{\lambda}\bigcap
U_{\iota}\not=\emptyset, \lambda, \iota\in\Lambda\}$ \vskip 2mm

\no where $\kappa_{\lambda\iota}$ is the number of non-homotopic
loops in formed between $U_{\lambda}$ and $U_{\iota}$. Then we get a
combinatorial manifold $\widetilde{M}_G(n)$ underlying a countable
graph $G$.

Define a combinatorial Euclidean space $\mathscr{E}_{G'}(n,
\lambda\in\Lambda)$ of spaces ${\bf R}^{n}$ by\vskip 3mm

\hskip 5mm $V(G')= \{\varphi_{\lambda}(U_{\lambda})|
\lambda\in\Lambda\},$\vskip 2mm

\hskip 5mm $E(G')=\{\
(\varphi_{\lambda}(U_{\lambda}),\varphi_{\iota}(U_{\iota}))_i,1\leq
i\leq\kappa'_{\lambda\iota}+1 | \
\varphi_{\lambda}(U_{\lambda})\bigcap
\varphi_{\iota}(U_{\iota})\not=\emptyset, \lambda,
\iota\in\Lambda\},$ \vskip 2mm

\no where $\kappa'_{\lambda\iota}$ is the number of non-homotopic
loops in formed between $\varphi_{\lambda}(U_{\lambda})$ and
$\varphi_{\iota}(U_{\iota})$. Notice that
$\varphi_{\lambda}(U_{\lambda})\bigcap
\varphi_{\iota}(U_{\iota})\not=\emptyset$ if and only if
$U_{\lambda}\bigcap U_{\iota}\not=\emptyset$ and
$\kappa_{\lambda\iota}=\kappa'_{\lambda\iota}$ for $\lambda,
\iota\in\Lambda$. We know that $G \ \cong \ G'$ by definition.

Now we prove that $\widetilde{M}_G(n)$ is homeomorphic to
$\mathscr{E}_{G'}(n, \lambda\in\Lambda)$. By assumption, $M^n$ is
locally compact. Whence, there exists a partition of unity
$c_{\lambda}: U_{\lambda}\rightarrow{\bf R}^n$, $\lambda\in\Lambda$
on the atlas $\mathscr{A}[M^n]$. Let $A_{\lambda}={\rm
supp}(\varphi_{\lambda})$. Define functions $h_{\lambda}:
M^n\rightarrow{\bf R}^n$ and ${\bf H}:
M^n\rightarrow\mathscr{E}_{G'}(n)$ by

\[
h_{\lambda}(x)=\left\{\begin{array}{lr}
c_{\lambda}(x)\varphi_{\lambda}(x) & {\rm if} \ x\in U_{\lambda},\\
{\bf 0}=(0,\cdots,0) & {\rm if} \ x\in U_{\lambda}-A_{\lambda}.
\end{array}
\right.
\]

\no and

$${\bf H}=\sum\limits_{\lambda\in\Lambda}\varphi_{\lambda}c_{\lambda}, \ \ \ {\rm
and} \ \ \ {\bf
J}=\sum\limits_{\lambda\in\Lambda}c_{\lambda}^{-1}\varphi_{\lambda}^{-1}.$$

\no Then $h_{\lambda}$, ${\bf H}$ and ${\bf J}$ all are continuous
by the continuity of $\varphi_{\lambda}$ and $c_{\lambda}$ for
$\forall\lambda\in\Lambda$ on $M^n$. Notice that
$c_{\lambda}^{-1}\varphi_{\lambda}^{-1}\varphi_{\lambda}c_{\lambda}=$the
unity function on $M^n$. We get that ${\bf J}={\bf H}^{-1}$, i.e.,
${\bf H}$ is a homeomorphism from $M^n$ to $\mathscr{E}_{G'}(n,
\lambda\in\Lambda)$. \hfill$\Box$\vskip 1mm

According to Theorem $2.1$, a combinatorial Euclidean space
homeomorphic to a $n$-manifold $M^n$ can be denoted by
$\mathscr{E}_{G[M^n]}(n,\ \mu\in\Lambda)$. We classify such
combinatorial Euclidean spaces $\mathscr{E}_G(n, \mu\in\Lambda)$
into two classes by considering the intersection
$\varphi_{\mu}^{-1}({\bf R}^n)\cap\varphi_{\nu}^{-1}({\bf R}^)$ for
$\mu, \nu\in\Lambda$ following:\vskip 3mm

\no{\bf Class 1.} \ \ For $\forall\mu,\nu\in\Lambda$,
$\varphi_{\mu}^{-1}({\bf R}^n)\cap\varphi_{\nu}^{-1}({\bf
R}^)=\emptyset$ or homeomorphic to ${\bf R}^n$.\vskip 2mm

\no{\bf Class 2.} \ There are $\mu,\nu\in\Lambda$ such that
$\varphi_{\mu}^{-1}({\bf R}^n)\cap\varphi_{\nu}^{-1}({\bf
R}^)\not=\emptyset$ with more than $2$ arcwise connected
components.\vskip 2mm

We respectively discuss Classes $1$ and $2$ by dimensional graphs
following.

\vskip 6mm

\no{\bf $2.1$ \ Dimensional Graphs $\widetilde{M}^n[G]$}

\vskip 4mm

\no Let $B^n=\{(x_1,x_2,\cdots,x_n)|\ (x_1,x_2,\cdots,x_n)\in{\bf
R}^n,\ x_1^n+x_2^n+\cdots+x_n^n<1\}$, $S^n=\{(x_1,x_2,\cdots,x_n)|\
(x_1,x_2,\cdots,x_n)\in{\bf R}^n,\ x_1^n+x_2^n+\cdots+x_n^n=1\}$ be
a $n$-dimensional open ball or $n$-sphere for an integer $n\geq 1$,
respectively. The combinatorial Euclidean spaces
$\mathscr{E}_G(n,\lambda\in\Lambda)$ in Class $1$ is such a
topological space that each non-homotopic loops comes from the graph
$G$, which can be characterized by graphs in spaces, i.e., {\it
$n$-dimensional graphs} defined following.

\vskip 4mm

\no{\bf Definition $2.2$} \ {\it An $n$-dimensional graph
$\widetilde{M}^n[G]$ is a combinatorial Euclidean space
$\mathscr{E}_G(n)$ of ${\bf R}^{n}_{\mu},\mu\in\Lambda$ underlying a
combinatorial structure $G$ such that}\vskip 2mm

($1$) {\it $V(G)$ is discrete consisting of $B^n$, i.e., $\forall
v\in V(G)$ is an open ball $B_v^n$;}\vskip 1mm

($2$) {\it $\widetilde{M}^n[G]\setminus V(\widetilde{M}^n[G])$ is a
disjoint union of open subsets $e_1,e_2,\cdots,e_m$, each of which
is homeomorphic to an open ball $B^n$;}\vskip 1mm

($3$) {\it the boundary $\overline{e}_i-e_i$ of $e_i$ consists of
one or two $B^n$ and each pair $(\overline{e}_i,e_i)$ is
homeomorphic to the pair $(\overline{B}^n,B^n)$;}\vskip 1mm

($4$) {\it a subset $A\subset\widetilde{M}^n[G]$ is open if and only
if $A\cap\overline{e}_i$ is open for $1\leq i\leq m$.}

\vskip 3mm

A {\it topological graph $\mathscr{T}[G]$} of a graph $G$ is a
$1$-dimensional graph in a topological space $\mathscr{P}$. We
restate it in the following.

\vskip 4mm

\no{\bf Definition $2.3$} \ {\it A topological graph
$\mathscr{T}[G]$ is a pair $(X,X^0)$ of a Hausdorff space $X$ with
its a subset $X^0$ such that}\vskip 2mm

($1$) {\it $X^0$ is discrete, closed subspaces of $X$;}\vskip 1mm

($2$) {\it $X-X^0$ is a disjoint union of open subsets
$e_1,e_2,\cdots,e_m$, each of which is homeomorphic to an open
interval $(0,1)$;}\vskip 1mm

($3$) {\it the boundary $\overline{e}_i-e_i$ of $e_i$ consists of
one or two points. If $\overline{e}_i-e_i$ consists of two points,
then $(\overline{e}_i,e_i)$ is homeomorphic to the pair
$([0,1],(0,1))$; if $\overline{e}_i-e_i$ consists of one point, then
$(\overline{e}_i,e_i)$ is homeomorphic to the pair
$(S^1,S^1-\{1\})$;}\vskip 1mm

($4$) {\it a subset $A\subset\widetilde{T}[G]$ is open if and only
if $A\cap\overline{e}_i$ is open for $1\leq i\leq m$.}

\vskip 3mm

Observation shows that there is a natural relation between a
$n$-dimensional graph $\widetilde{M}^n[G]$ with that of its a
topological graph $\mathscr{T}_0[G]=(X_0,X^0_0)$ which is
constructed from $\widetilde{M}^n[G]_1$ by:\vskip 3mm

($1$) Let $X^0_0=\{ \ O_v \ | \ O_v \ {\rm is \ the \ center \ of} \
B^n_v \ {\rm for} \ v\in V(G)\}$;\vskip 1mm

($2$) For $\forall uv\in E(G)$, let $uv$ be a line segment $e_{uv}:
[0,1]\rightarrow B^n_u\cup B^n_v$ such that $e_{uv}(0)=O_u$ and
$e_{uv}(1)=O_v$.\vskip 2mm

Then, the following result shows that a $n$-dimensional graph
$\widetilde{M}^n[G]$ is in fact a blown up of a topological graph
$\mathscr{T}[G]$ to dimensional $n$.

\vskip 4mm

\no{\bf Theorem $2.4$} \ {\it For any integer $n\geq 1$,
$\mathscr{T}_0[G]$ is a deformation retract of
$\widetilde{M}^n[G]$.}

\vskip 3mm

{\it Proof} \ If $n=1$, then $\widetilde{M}^n[G]=\mathscr{T}_0[G]$
is itself a topological graph. So we assume $n\geq 2$. Define a
mapping $f: \widetilde{M}^n[G]\times I\rightarrow\widetilde{M}^n[G]$
by

$$f(\overline{x},t)=(1-t)\overline{x}+t\overline{x}_0,$$

\no for $\forall\overline{x}\in\widetilde{M}^n[G]_1, t\in I$, where
$\overline{x}_0=O_v$ if $\overline{x}\in B^n_v$, and
$\overline{x}_0=p(\overline{x})$ if $\overline{x}\in e_i$, where $p:
uv\rightarrow e_{uv}$ a projection for $1\leq i\leq m$, such as
those shown in Fig.$2.1$. \vskip 3mm

\includegraphics[bb=0 0 180 200]{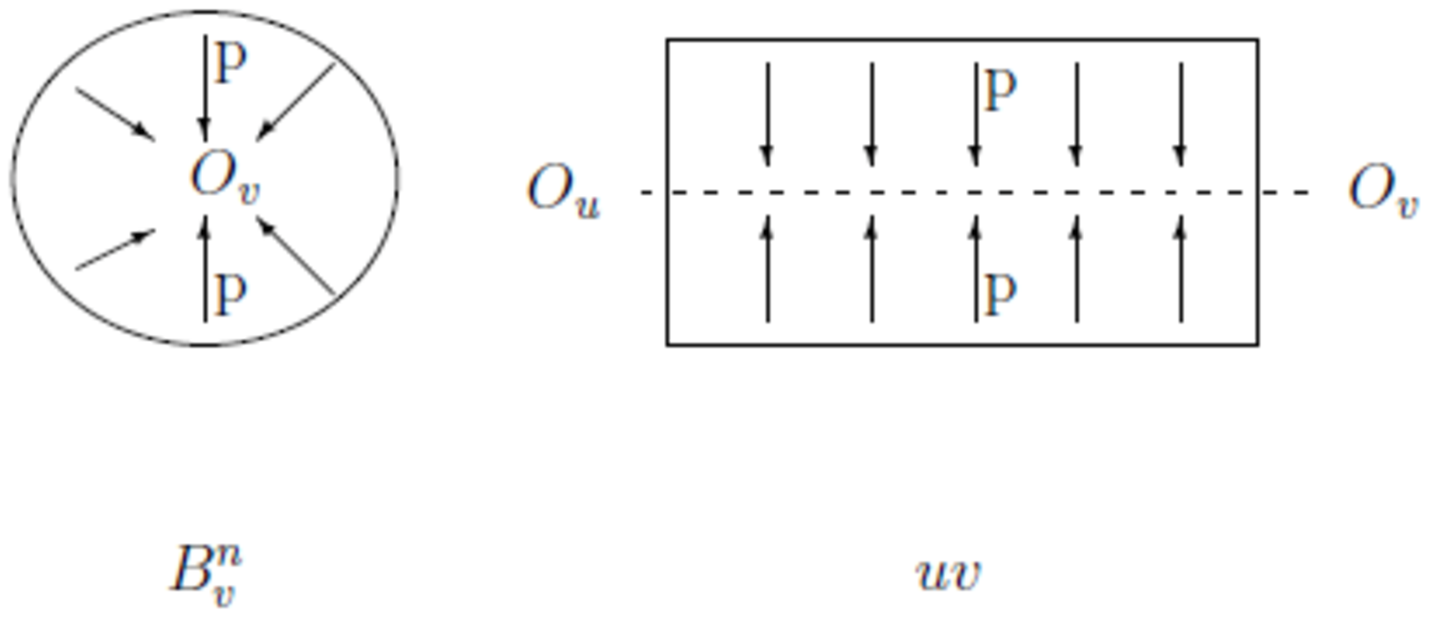}

\vskip 3mm

\c{\bf Fig.$2.1$}\vskip 2mm

\no Then we know that $f$ is continuous by definition and for
$\forall\overline{x}\in\widetilde{M}^n[G]$,

\vskip 3mm

\hskip 10mm$f(\overline{x},0)=\overline{x},$\vskip 2mm

\hskip 10mm$f(\overline{x},1)=p(\overline{x}),$\vskip 2mm

\no and $f(\overline{x},t)=\overline{x}$ for $\forall\overline{x}\in
\mathscr{T}_0[G]$ and $t\in I$. Therefore, $\mathscr{T}_0[G]$ is a
deformation retract of $\widetilde{M}^n[G]$ by
definition.\hfill$\Box$\vskip 1mm

Notice that the inclusion mapping $i:
\mathscr{T}_0[G]\rightarrow\widetilde{M}^n[G]$ is an isomorphism
between groups $\pi(\mathscr{T}_0[G],v_0)$ and
$\pi(\widetilde{M}^n[G]_1,v_0)$. We get a conclusion by Theorem
$2.4$ following.

\vskip 4mm

\no{\bf Corollary $2.5$} \ {\it Let $\widetilde{M}^n[G]$ be a
$n$-dimensional graph. Then for
$v_0\in\mathscr{T}_0[G]$,}

$$\pi(\widetilde{M}^n[G],v_0)\cong\pi(\mathscr{T}_0[G],v_0).$$

\vskip 3mm

We have known the structure of fundamental group of a topological
graph $\mathscr{T}[G]$ in [13]. Whence, we can characterize the
fundamental group of a $n$-dimensional graph $\widetilde{M}^n[G]$ by
applying Corollary $2.5$ following.

\vskip 4mm

\no{\bf Theorem $2.6$} \ {\it Let $T_{span}$ be a spanning tree in
the topological graph $\mathscr{T}_0[G]$, $\{e_{\lambda}: \
\lambda\in\Lambda\}$ the set of edges of $\mathscr{T}_0[G]$ not in
$T_{span}$ and
$\alpha_{\lambda}=A_{\lambda}e_{\lambda}B_{\lambda}\in\pi(\mathscr{T}_0[G],v_0)$
a loop associated with $e_{\lambda}=a_{\lambda}b_{\lambda}$ for
$\forall\lambda\in\Lambda$, where $v_0\in \mathscr{T}_0[G]$ and
$A_{\lambda}$, $B_{\lambda}$ are unique paths from $v_0$ to
$a_{\lambda}$ or from $b_{\lambda}$ to $v_0$ in $T_{span}$.
Then}

$$\pi(\mathscr{T}_0[G],v_0)=\left<\alpha_{\lambda}|\lambda\in\Lambda\right>.$$

An {\it $n$-dimensional tree} is a $n$-dimensional graph
$\widetilde{M}^n[G]$ with a tree $G$, denoted by
$\widetilde{M}^n[T]$. Applying Theorem $2.4$, we know the next
result.

\vskip 4mm

\no{\bf Theorem $2.7$} \ {\it An $n$-dimensional tree
$\widetilde{M}^n[T]$ is contractible.}

\vskip 3mm

{\it Proof} \ By Theorem $2.4$, we know that there is a continuous
mapping $f: \widetilde{M}^n[T]\times I\rightarrow\widetilde{M}^n[T]$
such that $\mathscr{T}_0[T]$ is a deformation retract of
$\widetilde{M}^n[T]$, i.e., for
$\forall\overline{x}\in\widetilde{M}^n[T]$,

\vskip 4mm

\hskip 10mm$f(\overline{x},0)=\overline{x},$\vskip 2mm

\hskip 10mm$f(\overline{x},1)=p(\overline{x}),$\vskip 3mm

\no and $f(\overline{x},t)=\overline{x}$ for $\forall\overline{x}\in
\mathscr{T}_0[T]$ and $t\in I$. Notice that $\mathscr{T}_0[T]$ is
contractible ([13]), we have a continuous mapping $g:
\mathscr{T}_0[T]\times I\rightarrow\mathscr{T}_0[T]$ such that
$\{v_0\}$ is a deformation retract for $\forall
v_0\in\mathscr{T}[T]$. Whence, the composition mapping $g\circ
f:\widetilde{M}^n[T]\times I\rightarrow\widetilde{M}^n[T]$ is
continuous such that for $\forall\overline{x}\in\widetilde{M}^n[T]$,

\vskip 4mm

\hskip 10mm$g\circ f(\overline{x},0)=\overline{x},$\vskip 2mm

\hskip 10mm$g\circ f(\overline{x},1)=p(\overline{x}),$\vskip 3mm

\no and $g\circ f(v_0,t)=v_0$ for $\forall t\in I$, i.e., $\{v_0\}$
is a deformation retract of $\widetilde{M}^n[T]_1$. This completes
the proof.\hfill$\Box$

\vskip 5mm

\no{\bf $2.2$ \ Labeled Dimensional Graphs
$\widetilde{M}^n[G^L]$}\vskip 3mm

\no In Class $2$, non-homotopic loops come from both
$\varphi_{\mu}^{-1}({\bf R}^n)\cap\varphi_{\nu}^{-1}({\bf R})$ for
$\mu,\nu\in\Lambda$ and the combinatorial structure $G$ of
$\mathscr{E}_G(n,\lambda\in\Lambda)$, which enables us to construct
the labeled dimensional graph $\widetilde{M}^n[G^L]$ following.

\vskip 4mm

\no{\bf Definition $2.8$} \ {\it Let $\mathscr{A}[M^n]=\{\
(U_{\lambda};\varphi_{\lambda})\ | \ \lambda\in\Lambda\}$ be an
atlas of an orientable manifold $M^n$. A labeled dimensional graph
$\widetilde{M}^n[G^L]$ of a combinatorial Euclidean space
$\mathscr{E}_G(n,\lambda\in\Lambda)$ is such a $n$-dimensional graph
$\widetilde{M}^n[G]$ with labeling $L:\
(U_{\mu},U_{\nu})\rightarrow\kappa_{\mu\nu}+1$ for
$\forall(U_{\mu},U_{\nu})\in E(G)$, where $\kappa_{\mu\nu}$ is the
number of non-homotopic loops between $U_{\mu}$ and $U_{\nu}$ for
$\mu,\nu\in\Lambda$.}\vskip 3mm

Notice that such an edge labeled graph $G^L$ can be simply viewed as
a graph $H_G$ with multiple edges by defining.\vskip 3mm

$V(H_G)=V(G^L)$;\vskip 2mm

$E(H_G)=\{(U_{\mu},U_{\nu})_i,\ 1\leq i\leq \kappa_{\mu\nu}+1\ |\
{\rm for} \ \forall(U_{\mu},U_{\nu})\in E(G),\
\mu,\nu\in\Lambda\}$.\vskip 3mm

\no{\bf Convention $2.9$} \ {\it An edge labeled graph $G^L$ with an
integer labeling $L(e)\geq 1$ for $\forall e\in E(G)$ is viewed as a
multiple graph $H_G$, i.e., $G^L=H_G$ throughout this paper.}\vskip
2mm

By this convention, such labeled graphs can be used essentially to
characterize non-homotopic loops of manifolds. For example, the
torus is a combinatorial Euclidean space
$\mathscr{E}_G(2,\lambda=1,2)$ underlying a dipole $D_{0,3,0}$,
i.e., consists of two Euclidean space $U_1,\ U_2\cong{\bf R}^2$ with
$\kappa_{ij}=2$ non-homotopic loops and its labeled graph on
$P_2^L=D_{0,3,0}$, such as those shown in Fig.$2.2$.

\vskip 3mm

\includegraphics[bb=0 0 180 180]{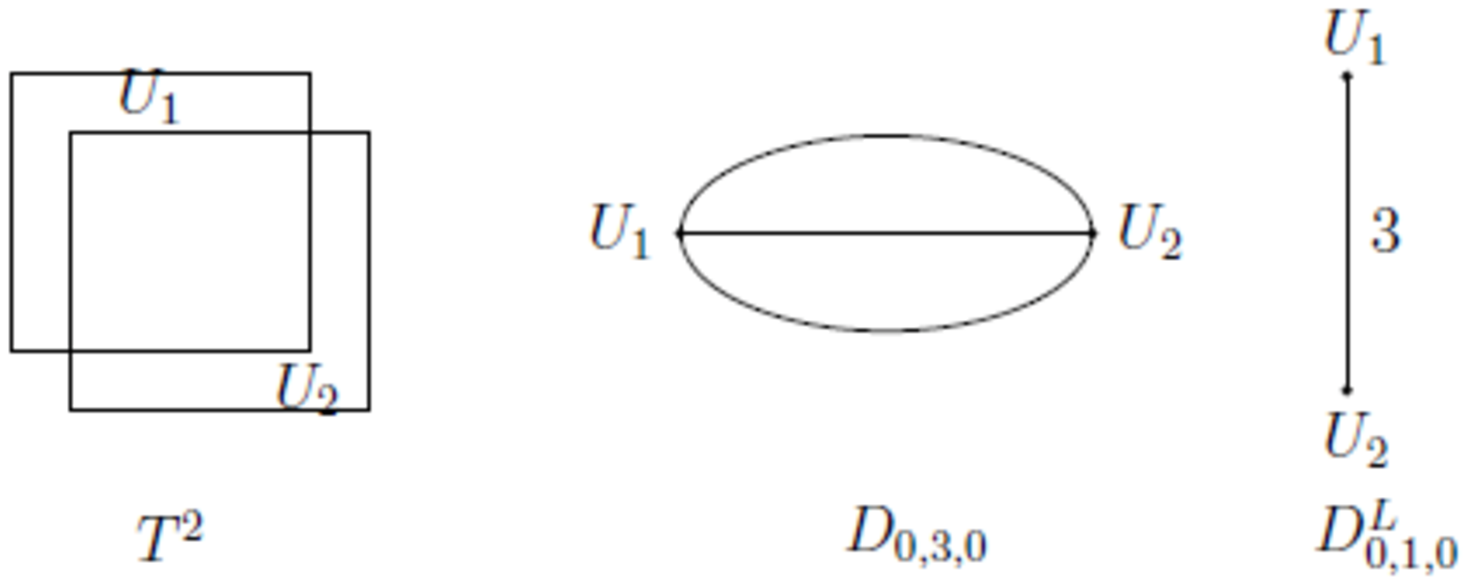}

\vskip 3mm

\c{\bf Fig.$2.2$}\vskip 2mm

Convention $2.9$ immediately implies the next result.

\vskip 4mm

\no{\bf Theorem $2.10$} \ {\it The number $\varpi(G^L)$ of cycles
basis of an edge labeled graph $G^L$ with an integer labeling
$L(e)\geq 1$ for $\forall e\in E(G)$ is\vskip 3mm

\c{$\varpi(G^L)=\varpi(G)+\sum\limits_{L(e)\geq 2, e\in
E(G)}(L(e)-1).$}\vskip 2mm

\no Whence, $\varpi(G^L)=0$ if and only if $\varpi(G)=0$ and
$L(e)=1$ for $\forall e\in E(G)$.}\vskip 3mm

Now denoted by $\mathscr{N}_C(\mathscr{E})$ and $\mathscr{N}_C(G^L)$
the sets of non-homotopic loops in the combinatorial Euclidean space
$\mathscr{E}_G(n,\lambda\in\Lambda)$ and cycles in its labeled
dimensional graph $\widetilde{M}^n[G^L]$, respectively. Then we know
the following interesting result.

\vskip 4mm

\no{\bf Theorem $2.11$} \ {\it There is a bijection $\vartheta:\
\mathscr{N}_C(\mathscr{E})\rightarrow\mathscr{N}_C(G^L)$, i.e.,}

$$\pi(\mathscr{E}_G(n,\lambda\in\Lambda)) \ \cong \ \pi(\widetilde{M}^n[G^L]).$$

\vskip 2mm

{\it Proof} \ We only need to prove that a loop
$L\in\mathscr{N}_C(\mathscr{E})$ if and only if there is a cycle
$C_L\in\mathscr{N}_C(G^L)$. The proof is divided into two cases
following. \vskip 3mm

\no{\bf Case 1.} \ $L$ comes from the combinatorial structure
$G$.\vskip 2mm

According to Theorems $2.1$ and $2.4$, we know that the underlying
graphs of $\mathscr{E}_G(n,\lambda\in\Lambda)$ and
$\widetilde{M}^n[G^L]$ are isomorphic. Whence there exists a cycle
$C_L\in\mathscr{T}_0[G]$, i.e., $C_L\in\mathscr{N}_C(G^L)$
correspondent to $L$ and verse via, for a cycle
$C\in\mathscr{T}_0[G]$, there also exists a loop
$L_C\in\mathscr{N}_C(\mathscr{E})$ for a cycle $C$. Whence, such a
mapping $\vartheta:\ L\rightarrow C_L$ is a bijection.\vskip 3mm

\no{\bf Case 2.} \ $L$ comes from $U_{\mu}\cap U_{\nu}$ for two
indexes $\mu,\nu\in\Lambda$.\vskip 2mm

Assume there are $\kappa_{\mu\nu}$ non-homotopic loos in $U\cap
U_{\nu}$. Not loss of generality, let $L$ be the {\it i}th loop. By
Definition $2.8$ and convention $2.9$, we have a cycle $C_L$
consisted of multiple edges $(U_{\mu},U_{\nu})_i$ with
$(U_{\mu},U_{\nu})_{i+1}$ in the graph $G^L$. Then $\vartheta:\
L\rightarrow C_L$ is a bijection by definition.

Combining the discussion of Cases $1$ and $2$, we get a
bijection

$$\vartheta:\
\mathscr{N}_C(\mathscr{E})\rightarrow\mathscr{N}_C(H_G).$$ Notice
that $\pi(\mathscr{E}_G(n))=\left<\mathscr{N}_C(\mathscr{E})\right>$
and
$\pi(\widetilde{M}^n[G^L])=\pi(\widetilde{M}^n[H_G])=\left<\mathscr{N}_C(G^L)\right>$.
The bijection $\vartheta:\
\mathscr{N}_C(\mathscr{E})\rightarrow\mathscr{N}_C(G^L)$ naturally
induces an isomorphism $\vartheta*:\
\pi(\mathscr{E}_G(n))\rightarrow\pi(\widetilde{M}^n[G^L])$ by
defining $\vartheta*(L_1+L_2)=\vartheta(L_1)+\vartheta(L_2)$ for
$L_1,\ L_2\in\mathscr{N}_C(\mathscr{E})$ in the field
$\mathbb{Z}_2$. Hence, we conclude that\vskip 4mm

\c{$\hskip 50mm\pi(\mathscr{E}_G(n)) \ \cong \
\pi(\widetilde{M}^n[G^L]).\hskip 50mm\Box$}\vskip 3mm

Combining Theorems $2.1$, $2.4$ and $2.11$, we get an important
result on $n$-manifold with finite non-homotopic loops following.

\vskip 4mm

\no{\bf Theorem $2.12$} \ {\it For an $n$-manifold $M^n$ with finite
non-homotopic loops, there always exists a labeled $n$-dimensional
graph $\widetilde{M}^n[G^L]$ such that}

$$\pi(M^n)\cong \pi(\widetilde{M}^n[G^L]).$$\vskip 2mm

{\it Proof} \ According to Theorem $2.1$, there is a combinatorial
Euclidean space $\mathscr{E}_G(n,\lambda\in\Lambda)$. Applying
Theorems $2.4$ and $2.11$, we get a labeled $n$-dimensional graph
$\widetilde{M}^n[G^L]$ such that

$$\pi(\mathscr{E}_G(n,\lambda\in\Lambda))\cong \pi(\widetilde{M}^n[G^L]),$$

\no where $G^L$ is defined in Definition $2.8$. Whence, we know that

$$\pi(M^n)\cong\pi(\mathscr{E}_G(n,\lambda\in\Lambda))\cong\pi(\widetilde{M}^n[G^L]).$$

\no This completes the proof. \hfill$\Box$\vskip 2mm

Applying Theorems $2.4$ and $2.12$, we get more useful conclusions
following.

\vskip 4mm

\no{\bf Corollary $2.13$} \ {\it The number of non-homotopic loops
in an $n$-manifold $M^n$ with finite non-homotopic loops is equal to
the dimension of cycle space of its
$\mathscr{T}_0[H_{G[M^n]}]$.}\vskip 3mm

\no{\bf Corollary $2.14$} \ {\it For an integer $n\geq 2$, a compact
$n$-manifold $M^n$ is homotopy sphere if and only if
$\widetilde{M}^n[G^L[M^n]]$ is a finite $n$-dimensional tree with
labeling $L:\ e\rightarrow 1$ for $\forall e\in E(G^L[M^n])$.}\vskip
3mm

{\it Proof} \ By definition, $M^n$ is compact if and only if $|G^L|$
is finite. Applying Theorems $2.10$ and $2.12$, we know that
$\pi(M^n)$ is trivial if and only if $\pi(\widetilde{M}^n[G^L])$ is
trivial, i.e., each labeling of edges is $1$. Now by Corollary
$2.5$, there must be $\pi(\mathscr{T}_0[H_G],v_0)=\{v_0\}$ for
$v_0\in\mathscr{T}_0[H_G]\cap M^n$. But this can happens only if
$\mathscr{T}_0[H_G]$ is a finite tree by Theorem $2.6$. Whence,
$\widetilde{M}^n[G^L[M^n]]$ is an $n$-dimensional tree with a
labeling $L(e)=1$ for $\forall e\in E(G^L[M^n])$.\hfill$\Box$\vskip
2mm

Corollary $2.14$ enables us to obtain an unified proof for the
generalized Poincar\'{e} conjecture, i.e., {\it any homotopy
$n$-sphere is homeomorphic to $S^n$} following.

\vskip 4mm

\no{\bf[Proof of Theorem $1.4$]}\vskip 3mm

\no For any integer $n\geq 1$, by Corollary $2.14$ an $n$-manifold
$M^n$ is homotopy $n$-sphere if and only if
$\widetilde{M}[G^L[M^n]]$ defined in Definition $2.8$ is a finite
labeled $n$-dimensional tree with $L(e)=1$ for $\forall e\in E(G)$,
i.e., a finite $n$-dimensional tree. Applying Theorem $2.7$, such an
$n$-dimensional tree $\widetilde{M}[G^L[M^n]]$ can deformation
retract to a point $v_0\in M^n$. Whence, $M^n$ is homeomorphic to an
$n$-sphere $S^n$.\hfill$\Box$\vskip 1mm

\vskip 8mm

\no{\bf\S$3.$ \ Listing $n$-Manifolds by Labeled Graphs}\vskip 5mm

\no{\bf Theorem $3.1$} \ {\it Let $\mathscr{A}[M^n]=\{\
(U_{\lambda};\varphi_{\lambda})\ | \ \lambda\in\Lambda\}$ be a atlas
of a locally compact $n$-manifold $M^n$. Then the dimensional graph
$\widetilde{M}^n[G^L_{|\Lambda}|]$ of $M^n$ is a topological
invariant on $|\Lambda|$, i.e., if
$\widetilde{M}^n[H^{L_1}_{|\Lambda|}]$ and
$\widetilde{M}^n[G^{L_2}_{|\Lambda|}]$ are two labeled
$n$-dimensional graphs of $M^n$, then there exists a
self-homeomorphism $h: M^n\rightarrow M^n$ such that $h:
H^{L_1}_{|\Lambda|}\rightarrow G^{L_2}_{|\Lambda|}$ naturally
induces an isomorphism of graph.}\vskip 3mm

{\it Proof} \ Let\vskip 3mm

\c{$\mathscr{A}^1_{|\Lambda|}[M^n]=\{\
(U_{\lambda};\varphi_{\lambda})\ | \ \lambda\in\Lambda_1\}$}\vskip
2mm

\no and\vskip 3mm

\c{$\mathscr{A}^2_{|\Lambda|}[M^n]=\{\ (V_{\lambda};\phi_{\lambda})\
| \ \lambda\in\Lambda_2\}$}\vskip 2mm

\no be two minimum atlases of a locally compact $n$-manifold $M^n$
with labeled graphs $H^{L_1}_{|\Lambda|}$ and $G^{L_2}_{|\Lambda|}$,
respectively, where $\Lambda_1=\Lambda_2=\{1,2,3,\cdots,k,\cdots\}$
are countable index sets. Notice that $\varphi_{\lambda}:
U_{\lambda}\rightarrow{\bf R}^n$ and $\phi_{\lambda}:
V_{\iota}\rightarrow{\bf R}^n$ are homeomorphisms for $\forall
\lambda\in\Lambda_1$ and $\iota\in\Lambda_2$. So there is a
homeomorphism $\tau_{\lambda}:
\varphi_{\lambda}(U_{\lambda})\rightarrow\phi_{\lambda}(V_{\lambda})$
for $\lambda\in\Lambda_1=\Lambda_2$. Now define

$$h_{\lambda}=\phi_{\lambda}^{-1}\tau_{\lambda}\varphi_{\lambda}: x\rightarrow
\phi_{\lambda}^{-1}(\tau_{\lambda}(\varphi_{\lambda}(x))) \ \ {\rm
for} \ \ x\in U_{\lambda}.$$

\no Then $h_{\lambda}$ with its inverse
$h^{-1}=\varphi_{\lambda}^{-1}\tau_{\lambda}^{-1}\phi_{\lambda}$ is
continuous on $M^n$. By the compactness of $M^n$ there exists a
partition of unity $c_{\lambda}: U_{\lambda}\rightarrow{\bf R}^n$,
$\lambda\in\Lambda$ on the atlas $\mathscr{A}^1_{|\Lambda|}[M^n]$.
Define a function $h: M^n\rightarrow M^n$ by

$$h =\sum\limits_{\lambda\in\Lambda} \phi_{\lambda}^{-1}\tau_{\lambda}\varphi_{\lambda}c_{\lambda},$$

\no where $A_{\lambda}={\rm supp}(\varphi_{\lambda})$ and\vskip 3mm

\[
\phi_{\lambda}^{-1}\tau_{\lambda}\varphi_{\lambda}c_{\lambda}(x)=\left\{\begin{array}{lr}
c_{\lambda}(x)\phi_{\lambda}^{-1}\tau_{\lambda}\varphi_{\lambda}(x) & {\rm if} \ x\in U_{\lambda},\\
{\bf 0}=(0,\cdots,0) & {\rm if} \ x\in U_{\lambda}-A_{\lambda}.
\end{array}
\right.
\]\vskip 2mm

\no Then $h: M^n\rightarrow M^n$ is a homeomorphism with
$h(U_{\lambda})=V_{\lambda}$ for $\lambda\in\Lambda_1=\Lambda_2$.
Hence, $h:\ V(H_{|\Lambda|}^L)\rightarrow V(G^L_{|\Lambda|})$ is a
bijection by definition.

Now if there are $\kappa_{\mu\nu}$ non-homotopic loops between
$U_{\mu}$ and $U_{\nu}$, then there are must be $\kappa_{\mu\nu}$
non-homotopic loops between $V_{\mu}$ and $V_{\nu}$ and vice via by
the homeomorphic property. Therefore, $L_{1}:
(U_{\mu},U_{\nu})\rightarrow\kappa_{\mu\nu}+1$ in $H^L_{|\Lambda|}$
if and only if $L_2: (V_{\mu},V_{\nu})\rightarrow\kappa_{\mu\nu}+1$
in $G^{L_2}_{|\Lambda|}$, i.e.,
$h(U_{\mu},U_{\nu})=(h(U_{\mu}),h(V_{\nu}))$ with
$L_1(U_{\mu},U_{\nu})=L_2(h(U_{\mu}),h(U_{\nu}))$. By definition,
two labeled graphs $G^L_1$ and $G^L_2$ with labeling mappings $L_1$
and $L_2$ are said to be isomorphic if there is an isomorphism
$\varpi:\ G_1\rightarrow G_2$ with $\varpi L_1=L_2\varpi$. Whence,
$h:\ H^L_{|\Lambda|}\rightarrow G^L_{|\Lambda|}$ naturally induces
an isomorphism between labeled graphs $H^{L_1}_{|\Lambda|}$ and
$G^{L_2}_{|\Lambda|}$.\hfill$\Box$\vskip 2mm

We get a conclusion by Theorem $3.1$ following.

\vskip 4mm

\no{\bf Corollary $3.2$} \ {\it The labeled graph $G^L_{|\Lambda|}$
of a locally compact $n$-manifold $M^n$ is unique dependent on
$|\Lambda|$.}\vskip 3mm

For classifying $n$-manifolds by applying Theorem $3.1$, we
introduce the conception of minimum atlas following.

\vskip 4mm

\no{\bf Definition $3.3$} \ {\it An atlas\vskip 3mm

\c{$\mathscr{A}[M^n]=\{\ (U_{\lambda};\varphi_{\lambda})\ | \
\lambda\in\Lambda\}$}\vskip 2mm

\no of an $n$-manifold $M^n$ is minimal if there are no indexes
$\mu,\ \nu\in\Lambda$ and a continuous mapping $\varphi_{\mu\nu}$
with $\varphi_{\mu\nu}:\ U_{\mu}\cup U_{\nu}\rightarrow{\bf R}^n$
such that

\vskip 3mm

\c{$\mathscr{A}'=\{\ (U_{\lambda};\varphi_{\lambda}), (U_{\mu}\cup
U_{\nu},\varphi_{\mu\nu})\ | \
\lambda\in\Lambda\setminus\{\mu,\nu\}\}$}\vskip 2mm

\no is also an atlas of $M^n$.

An atlas $\mathscr{A}[M^n]$ of an $n$-manifold $M^n$ is minimum if
it has minimum cardinality among all of its minimal atlases. Denoted
such a minimal atlas by $\mathscr{A}_{min}[M^n]$ and its labeled
$n$-dimensional graph by $\widetilde{M}^n[G^L_{min}]$.}\vskip 3mm

The next result characterizes minimal atlases of an $n$-manifold.

\vskip 4mm

\no{\bf Theorem $3.4$} \ {\it Let \vskip 3mm

\c{$\mathscr{A}[M^n]=\{\ (U_{\lambda};\varphi_{\lambda})\ | \
\lambda\in\Lambda\}$}\vskip 2mm

\no be an atlas of a locally compact $n$-manifold $M^n$. Then}\vskip
3mm

($i$) {\it $\mathscr{A}[M^n]$ is minimal if and only if there are no
indexes $\mu,\ \nu\in\Lambda$ such that $U_{\mu}\cap U_{\nu}$ is
arcwise connected, i.e., $L(e)\geq 2$ for $\forall e\in
E(G^L_{min})$.}\vskip 1mm

($ii$) {\it $M^n$ is with finite non-homotopic loops if and only if
$\mathscr{A}_{min}[M^n]$ is finite.}

\vskip 3mm

{\it Proof} \ ($i$) \ If the result ($i$) is not true, then there
exist indexes $\mu,\ \nu\in\Lambda$ such that $U_{\mu}\cap U_{\nu}$
is arcwise connected. Assume\vskip 3mm

\c{$\varphi_{\mu}(U_{\mu}\cap U_{\nu})= S\subset{\bf R}^n \ \ {\rm
and} \ \ \varphi_{\nu}(U_{\mu}\cap U_{\nu})= T\subset{\bf
R}^n.$}\vskip 2mm

Notice that $S$ and $T$ are homeomorphic to ${\bf R}^n$ by
definition. We can always choose a continuous mapping $\tau:\
S\rightarrow T$, i.e., $\tau(S)=T$. Define\vskip 3mm

\c{$\varphi_{\mu}^{\tau}=\tau\varphi_{\mu}$.}\vskip 2mm

\no Then we get that $\varphi_{\mu}^{\tau}|_{U_{\mu}\cap
U_{\nu}}=\varphi_{\nu}|_{U_{\mu}\cap U_{\nu}}$. Whence, there is a
continuous mapping $\varphi_{\mu\nu}:\ U_{\mu}\cup
U_{\nu}\rightarrow{\bf R}^n$. Therefore,

\vskip 3mm

\c{$\mathscr{A}'=\{\ (U_{\mu}\cup U_{\nu}; \varphi_{\mu\nu}),
(U_{\lambda};\varphi_{\lambda})\ |\
\lambda\in\Lambda\setminus\{\mu,\nu\}\}$}\vskip 2mm

\no is also an atlas of $M^n$ but with $|\mathscr{A}'|=|\Lambda|-1$.
This contradicts to the minimality of $\Lambda$. So ($i$)
holds.\vskip 1mm

($ii$) \ If $\mathscr{A}_{min}[M^n]$ is finite, then $M^n$ is
obvious only with finite non-homotopic loops by the assumption of
its locally compactness. Now let

\vskip 3mm

\c{$\mathscr{A}_{min}[M^n]=\{\ (U_{\lambda};\varphi_{\lambda})\ | \
\lambda\in\Lambda\}$}\vskip 2mm

\no be a minimum atlas of a locally compact manifold $M^n$ with an
infinite index set $\Lambda$. By ($i$), there are no indexes $\mu,
\nu$ in $\Lambda$ such that $U_{\mu}\cap U_{\nu}$ is arcwise
connected. In other words, $U_{\mu}\cap U_{\nu}=\emptyset$ or with
more than $2$ arcwise components for $\forall\mu, \nu\in\Lambda$,
i.e., $L(e)\geq 2$ for $\forall e\in E(G^L_{min})[M^n]$. Applying
Theorem $2.10$, we know that the number

$$\varpi(G^L_{min}[M^n])=\varpi(G_{min}[M^n])+\sum\limits_{e\in E(G_{min}[M^n])}(L(e)-1)$$

\no of cycle basis of $G^L_{min}[M^n]$ is greater than any
sufficient larger number $N>0$, which contradicts the assumption
that $M^n$ is only with finite non-homotopic loops.
\hfill$\Box$\vskip 1mm

Combining Corollary $3.2$ with that of Theorem $3.4$, we get a
conclusion following, which enables us to list locally compact
$n$-manifolds with finite non-homotopic loops by labeled graphs.

\vskip 4mm

\no{\bf Corollary $3.5$} \ {\it If the minimum labeled graph
$G^L_{min}[M^n_1]$ of a locally compact $n$-manifold $M^n_1$ is not
isomorphic to the minimum labeled graph $G^L_{min}[M^n_2]$ of
$M_2^n$, then $M^n_1$ is not homeomorphic to $M^n_2$.}\vskip 3mm

Now by Theorem $3.4$, let\vskip 3mm

\c{$\mathscr{A}_{min}[M^n]=\{\ (U_{\lambda};\varphi_{\lambda})\ | \
\lambda\in\Lambda, |\Lambda|<+\infty\}$}\vskip 2mm

\no be a minimum atlas of locally compact $n$-manifolds $M^n$. Then
we can list $n$-manifolds following by Corollary $3.5$.

\vskip 3mm

{\bf (1) \ $|\Lambda|=1$}

\vskip 2mm

In this case, $\mathscr{A}_{min}[M^n]=\{(U;\varphi)\}$, i.e.,
$M^n={\bf R}^n$.\vskip 3mm

{\bf (2) \ $|\Lambda|=2$}\vskip 2mm

In this case,

$$\mathscr{A}_{min}[M^n]=\{(U_1;\varphi_1),\ (U_2;\varphi_2)\}$$

\no and $M^n$ is double covered, which can be classified by labeled
graphs $D^L_{0,1,0}$ shown in Fig.$3.1$,\vskip 3mm

\vskip 3mm

\c{\includegraphics[bb=0 0 180 60]{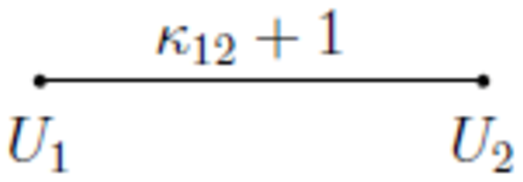}}

\vskip 3mm\vskip 3mm

\c{\bf Fig.$3.1$}\vskip 2mm

For example, if $n=2$, i.e., compact $2$-manifolds ${\bf S}$, then
$\kappa_{12}$ is the genus of ${\bf S}$ with $\kappa_{12}\geq 1$.
\vskip 3mm

{\bf (3) \ $|\Lambda|=3$}\vskip 2mm

In this case,

$$\mathscr{A}_{min}[M^n]=\{(U_1;\varphi_1),\ (U_2;\varphi_2),\
(U_3;\varphi_3)\}$$

\no and we can list such $n$-manifolds $M^n$ by labeled graphs in
Fig.$3.2$, where integers $k,l,s\geq 2$.\vskip 3mm

\vskip 3mm

\includegraphics[bb=0 0 180 120]{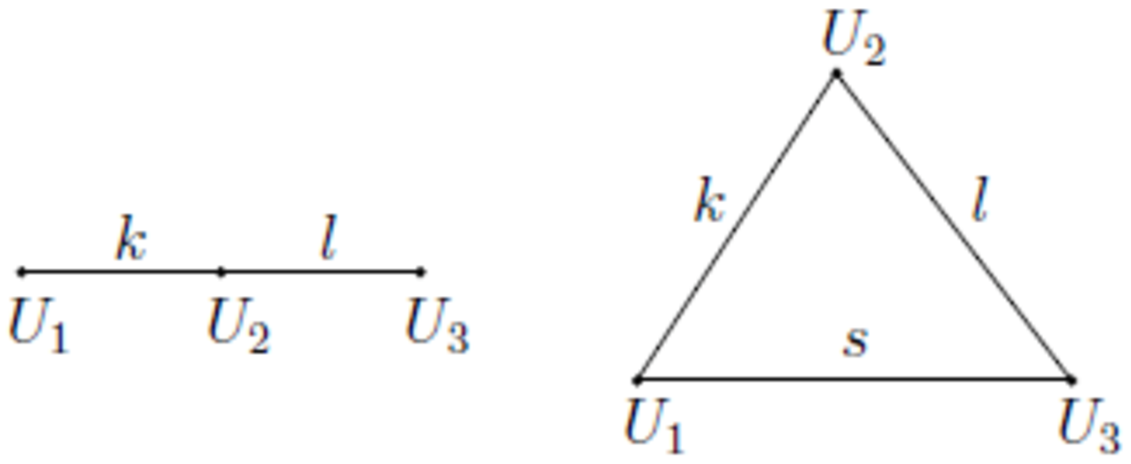}

\vskip 3mm\vskip 3mm

\c{\bf Fig.$3.2$}\vskip 3mm

{\bf (4) \ $|\Lambda|=4$}\vskip 2mm

In this case,

$$\mathscr{A}_{min}[M^n]=\{(U_1;\varphi_1),\ (U_2;\varphi_2),\
(U_3;\varphi_3),\ (U_4;\varphi_4)\}$$

\no and we can list such $n$-manifolds $M^n$ by labeled graphs in
Fig.$3.3$, where integers $k,l,r,q,s,t\geq 2$.\vskip 3mm

\vskip 3mm

\includegraphics[bb=0 0 180 315]{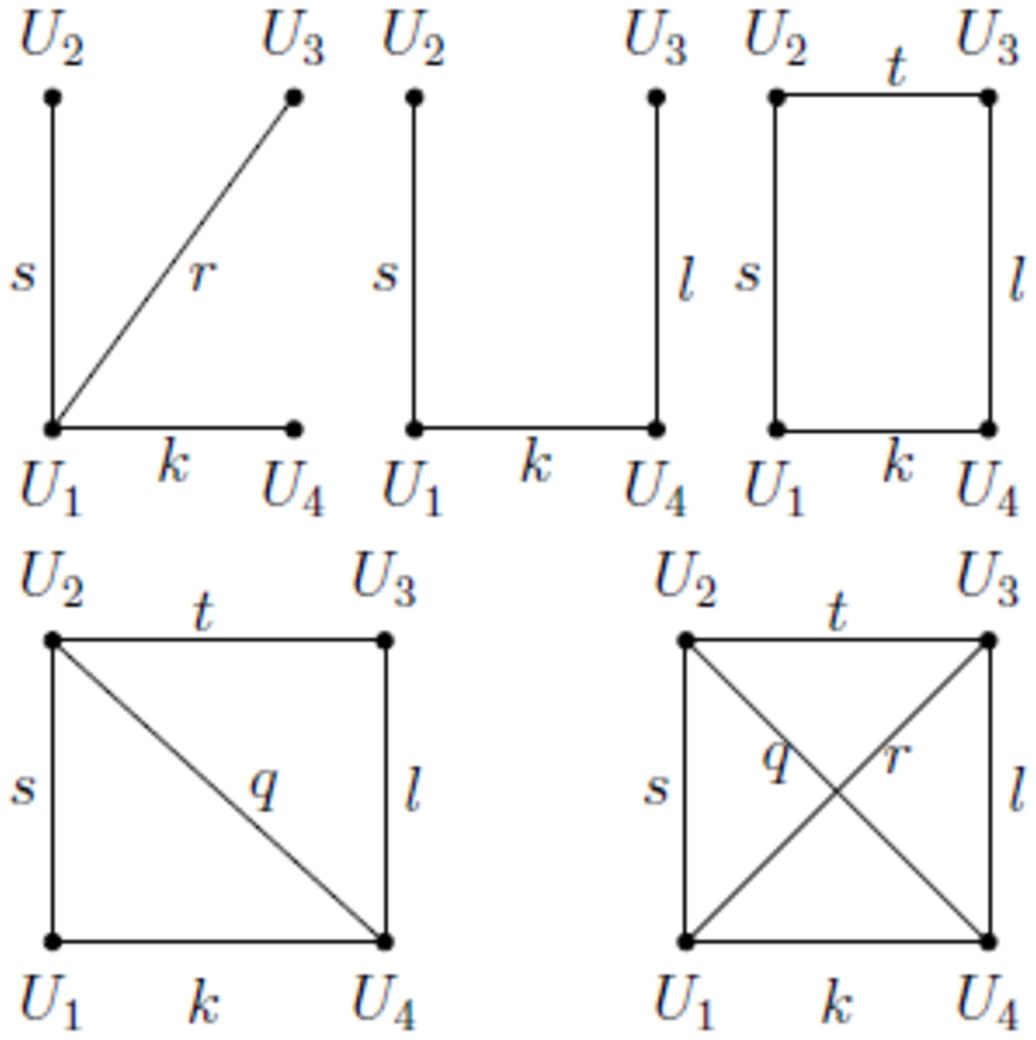}

\c{\bf Fig.$3.3$}\vskip 3mm

$\cdots\cdots\cdots\cdots$\vskip 3mm

{\bf (p) \ $|\Lambda|=p$}\vskip 2mm

In this case,

$$\mathscr{A}_{min}[M^n]=\{(U_i;\varphi_i), 1\leq i\leq p\}.$$

\no and there are $G_1,G_2,\cdots,G_{k(p)}$ such non-isomorphic
graphs known in graph theory, where

$$k(p)\sim
\frac{2^{\displaystyle\frac{p(p-1)}{2}}}{n!}$$

\no is the number of non-isomorphic graphs of order $p$ with

\vskip 0.2cm
\begin{center}
\begin{tabular}{c|ccccccccc}
$p$ & $1$ &  $2$ & $3$ & $4$ & $5$ & $6$ & $7$ & $8$ & $9$  \\
\hline
$k(p)$ & $1$ &  $1$ & $2$ & $6$ & $21$ & $112$ & $853$ & $11117$ & $261080$  \\
\end{tabular}
\end{center}
\vskip 0.2cm

\no for $p\leq 9$. Then we can list such $n$-manifolds $M^n$ by
labeled graphs following:

$$G^L_1,G^L_2,\cdots,G^L_{k(p)},\ \ \ L_i(e_{ij})\geq 2 \ \ {\rm for} \ \ \forall e_{i}\in E(G_i), 1\leq i\leq k(p).$$

\vskip 12mm

\no{\bf References}\vskip 3mm

\re{[1]}R.Abraham and J.E.Marsden, {\it Foundation of Mechanics}(2nd
edition), Addison-Wesley, Reading, Mass, 1978.

\re{[2]}M.H.Freedman, The topology of four-dimensional manifolds,
{\it J.Diff.Geom.}, 17 (1982), 357-453.

\re{[3]}S.Gadgil and H.Seshadri, Ricci flow and Perelman's proof of
the Poincar\'{e} conjecture, {\it Current Science}, Vol.91, No.10,
25 Nov. 2006, 1326-1334.

\re{[5]}J.L.Gross and T.W.Tucker, {\it Topological Graph Theory},
John Wiley \& Sons, 1987.

\re{[6]}L.F.Mao, {\it Automorphism Groups of Maps, Surfaces and
Smarandache Geometries}, American Research Press, 2005.

\re{[7]}L.F.Mao, {\it Smarandache Multi-Space Theory}, Hexis.
Phoenix, USA 2006.

\re{[8]}L.F.Mao, Geometrical theory on combinatorial manifolds, {\it
JP J.Geometry and Topology}, Vol.7, No.1(2007),65-114.

\re{[9]}L.F.Mao, Combinatorial speculation and combinatorial
conjecture for mathematics, {\it International J.Math. Combin.},
Vol.1(2007), 1-19.

\re{[10]}L.F.Mao, Combinatorial fields - an introduction, {\it
International J.Math. Combin.} Vol.3 (2009), 01-22.

\re{[11]}L.F.Mao, {\it Combinatorial Geometry with Applications to
Field Theory}, InfoQuest, USA, 2009.

\re{[12]}L.F.Mao, A combinatorial decomposition of Euclidean space
${\bf R}^n$ with contribution to visibility, {\it International
J.Math.Combin.}, Vol.1(2009).

\re{[13]}W.S.Massey, {\it Algebraic Topology: An Introduction},
Springer-Verlag, New York, etc.(1977).

\re{[14]}J.Milnor, Towards the Poincar\'{e} conjecture and the
classification of $3$-manifolds, {\it Notices AMS}, 50(2003),
1226-1233.

\re{[15]}G.Perelman, The entropy formula for the Ricci flow and its
geometric applications, {\it arXiv: math.DG/0211159v1}, 11 Nov.
2002.

\re{[16]}G.Perelman, Ricci flow with surgery on 3-manifolds, {\it
arXiv: math.DG/0303109}, 10 Mar 2003.

\re{[17]}G.Perelman, Finite extinction time for the solutions to the
Ricci flow on certain three-manifolds, {\it arXiv: math.DG/0307245},
17 Jul. 2003.

\re{[18]}S.Smale, Generalized Poincar\'{e} conjecture in dimensions
greater than four, {\it Annals of Math.}, 74(1961), 391-406.

\re{[19]}W.P.Thurston, {\it Three-Dimensional Geometry and
Topology}, Vol.1, ed. by Silvio Levy, Princeton Mathematical Series
35, Princeton University Press, 1997.

\end{document}